\documentclass[12pt]{article}

\usepackage{amsfonts,amsmath}

\def \beq {\begin{eqnarray}}
\def \eeq {\end{eqnarray}}
\def \beqn {\begin{eqnarray*}}
\def \eeqn {\end{eqnarray*}}

\newcommand{\halmos}{\rule{1ex}{1.4ex}}

\newcounter{for}[section]

\newtheorem{itlemma}{Lemma}[section]
\newtheorem{itproposition}[itlemma]{Proposition}
\newtheorem{theorem}[itlemma]{Theorem}
\newtheorem{itcorollary}[itlemma]{Corollary}
\newtheorem{itremark}[itlemma]{Remark}
\newtheorem{itremarks}[itlemma]{Remarks}
\newtheorem{itdefinition}[itlemma]{Definition}
\newtheorem{itexample}[itlemma]{Example}

\newenvironment{fact}{\begin{itfact}\rm}{\end{itfact}}
\newenvironment{claim}{\begin{itclaim}\rm}{\end{itclaim}}
\newenvironment{lemma}{\begin{itlemma}}{\end{itlemma}}
\newenvironment{remark}{\begin{itremark}\rm}{\end{itremark}}
\newenvironment{remarks}{\begin{itremarks} \rm}{\end{itremarks}}
\newenvironment{corollary}{\begin{itcorollary}}{\end{itcorollary}}
\newenvironment{proposition}{\begin{itproposition}}{\end{itproposition}}
\newenvironment{definition}{\begin{itdefinition}\rm}{\end{itdefinition}}
\newenvironment{example}{\begin{itexample}\rm}{\end{itexample}}
\newenvironment{proof}{\noindent {\em Proof}.\ \
}{\hspace*{\fill}$\halmos$\medskip}
\newcommand{\be}[1]{\addtocounter{for}{1} \begin{equation}\label{#1}}
\newcommand{\ee}{\end{equation}}

\newcommand{\bl}[1]{\begin{lemma}\label{#1}}
\newcommand{\br}[1]{\begin{remark}\label{#1}}
\newcommand{\brs}[1]{\begin{remarks}\label{#1}}
\newcommand{\bt}[1]{\begin{theorem}\label{#1}}
\newcommand{\bd}[1]{\begin{definition}\label{#1}}
\newcommand{\bp}[1]{\begin{proposition}\label{#1}}
\newcommand{\bc}[1]{\begin{corollary}\label{#1}}
\newcommand{\bfact}[1]{\begin{fact}\label{#1}}
\newcommand{\bex}[1]{\begin{example}\label{#1}}

\newcommand{\ec}{\end{corollary}}
\newcommand{\efact}{\end{fact}}
\newcommand{\eex}{\end{example}}
\newcommand{\el}{\end{lemma}}
\newcommand{\er}{\end{remark}}
\newcommand{\ers}{\end{remarks}}
\newcommand{\et}{\end{theorem}}
\newcommand{\ed}{\end{definition}}
\newcommand{\ep}{\end{proposition}}
\newcommand{\epr}{\end{proof}}
\newcommand{\bpr}{\begin{proof}}
\newcommand{\bcl}[1]{\begin{claim}\label{#1}}
\newcommand{\ecl}{\end{claim}}

\newcommand{\ecs}{\end{corollary}}
\newcommand{\eers}{\end{exercise}}
\newcommand{\eexs}{\end{example}}
\newcommand{\eems}{\end{example}}
\newcommand{\els}{\end{lemma}}
\newcommand{\eles}{\end{lemmaex}}
\newcommand{\ets}{\end{theorem}}
\newcommand{\eds}{\end{definition}}
\newcommand{\eps}{\end{proposition}}

\newcommand{\bi}{\begin{itemize}}
\newcommand{\ei}{\end{itemize}}
\newcommand{\ben}{\begin{enumerate}}
\newcommand{\een}{\end{enumerate}}

\def\vbar{\mathchoice{\vrule height6.3ptdepth-.5ptwidth.8pt\kern-.8pt}
   {\vrule height6.3ptdepth-.5ptwidth.8pt\kern-.8pt}
   {\vrule height4.1ptdepth-.35ptwidth.6pt\kern-.6pt}
   {\vrule height3.1ptdepth-.25ptwidth.5pt\kern-.5pt}}
\def\fudge{\mathchoice{}{}{\mkern.5mu}{\mkern.8mu}}
\def\bbc#1#2{{\rm \mkern#2mu\vbar\mkern-#2mu#1}}
\def\bbb#1{{\rm I\mkern-3.5mu #1}}
\def\bba#1#2{{\rm #1\mkern-#2mu\fudge #1}}
\def\bb#1{{\count4=`#1 \advance\count4by-64 \ifcase\count4\or\bba A{11.5}\or
   \bbb B\or\bbc C{5}\or\bbb D\or\bbb E\or\bbb F \or\bbc G{5}\or\bbb H\or
   \bbb I\or\bbc J{3}\or\bbb K\or\bbb L \or\bbb M\or\bbb N\or\bbc O{5} \or
   \bbb P\or\bbc Q{5}\or\bbb R\or\bbc S{4.2}\or\bba T{10.5}\or\bbc U{5}\or
   \bba V{12}\or\bba W{16.5}\or\bba X{11}\or\bba Y{11.7}\or\bba Z{7.5}\fi}}

\def \qed {{\hspace*{\fill}$\halmos$\medskip}}
\def \C {{\cal{C}}}

\def \s {\sigma}

\def \LL {{\cal{L}}}

\def \M {{\cal{M}}}
\def \AA {{\Re}}

\def \Pn {P_{\nu_{\rho}}(\tau>t)}
\def \Pe {P_{\eta}^*(\tau>t)}
\def \Pj {P_{\AA_{j}\eta}^*(\tau>t)}
\def \Pi {P_{\AA_i\eta}^*(\tau>t)}
\def \F {{\cal{F}}}

\newcommand{\ba}[1]{\addtocounter{for}{1} \begin{eqnarray}\label{#1}}
\newcommand{\ea}{\end{eqnarray}}

\def\sqr#1#2{{\vcenter{\vbox{\hrule height .#2pt
                             \hbox{\vrule width .#2pt height#1pt \kern#1pt
                                   \vrule width .#2pt}
                             \hrule height .#2pt}}}}

\def\pmb#1{\setbox0=\hbox{#1}%
   \kern-.025em\copy0\kern-\wd0
   \kern.05em\copy0\kern-\wd0
   \kern-.025em\raise.0433em\box0 }
\def\sqr#1#2{{\vcenter{\vbox{\hrule height.#2pt
     \hbox{\vrule width.#2pt height#1pt \kern#1pt
   \vrule width.#2pt}\hrule height.#2pt}}}}

\def\ZZ{{\mathbb Z}}   %
\def\A{{\mathcal A}}
\def\B{{\mathcal B}}
\def\RR{{\mathbb R}}   %
\def\NN{{\mathbb N}}   %
\def\PP{{\mathbb P}}   %

\def\d{\delta}
\def\D{\Delta}
\def\l{\lambda}
\def\e{\epsilon} 
\def\g{\gamma}
\def\L{\Lambda}
\def\nur{\nu_{\rho}}
\def\mur{\mu_{\rho}}
\def\G{\Gamma}

\def\v{\varphi}
\def\p{\partial}
\def\bs{\backslash}

\begin{document}

\title{Existence of quasi-stationary measures for asymmetric
attractive particle systems on $\ZZ^d$.}
\author{Amine Asselah \& Fabienne Castell\\
Laboratoire d'Analyse, Topologie et Probabilit\'es.  CNRS UMR 6632.\\
C.M.I., Universit\'e de Provence,\\
39 Rue Joliot-Curie, \\F-13453 Marseille cedex 13, France\\
asselah@gyptis.univ-mrs.fr \& castell@gyptis.univ-mrs.fr}
\date{}
\maketitle
\begin{abstract}
We show the existence of non-trivial quasi-stationary measures
for conservative attractive particle systems on $\ZZ^d$ conditioned
on avoiding an increasing local set $\A$. Moreover,
we exhibit a sequence of measures $\{\nu_n\}$, whose $\omega$-limit set
consists of quasi-stationary measures. For zero range processes, 
with stationary measure $\nur$, we prove
the existence of an $L^2(\nur)$ nonnegative eigenvector for the generator
with Dirichlet boundary on $\A$, after establishing
a priori bounds on the $\{\nu_n\}$.
\end{abstract}



{\em Running head}: Q-S measures for attractive systems.

\section{Introduction}
We consider the `processus des misanthropes', which includes
the asymmetric exclusion process and zero range processes. 
For concreteness,
let us describe here the dynamics of a zero range process. 
We denote the path of the process by $\{\eta_t,t\ge 0\}$ with
$\eta_t(i) \in \NN$ for $i\in \ZZ^d$.
At site $i$ and at time $t$, one of the $\eta_t(i)$ particles
jumps to site $j$ at rate $g(\eta_t(i))p(i,j)$ where 
\be{eq1.0}
g:\NN  \to [0,\infty)  \text{ is nondecreasing, with }
g(0)=0, \quad \sup_k\left(g(k+1)-g(k)\right)<\infty, 
\end{equation} 
and $p(.,.)$ is the transition kernel of a transient random walk.
Under assumptions that we make precise later, the informal
dynamics described above corresponds to a Feller process
with stationary product measures $\{\nu_{\rho},\rho>0\}$ (see \cite{andjel}).

Our motivation stems from statistical physics where such systems model gas
of charged particles in equilibrium under an electrical field. An interesting 
issue is the distribution of the occurrence time of density fluctuations
in equilibrium. Thus, let $\Lambda$ be a finite subset of $\ZZ^d$ and 
consider the event
\be{eq1.2}
\A=\{\eta:\ \frac{1}{|\Lambda|} \sum_{i\in \Lambda} \eta(i)>\rho'\}
\quad{\rm with}\quad \rho'>\rho.
\end{equation} 
Let $\tau$ be the first time a trajectory $\{\eta_t:t\ge 0\}$ enters $\A$.
As in \cite{ad,af}, we consider two complementary issues: 
\begin{itemize} 
\item[(i)] to estimate the tail of the distribution of $\tau$;
\item[(ii)] to characterize the law of $\eta_t$ at large time, conditioned on
$\{\tau>t\}$, when the initial configurations are drawn from $\nur$.
\end{itemize} 

We denote by $\LL$ the generator of our process, by $\{S_t,t\ge 0\}$
the associated semi-group, and by $P_{\mu}$ the law of the process
with initial probability $\mu$. 
For any probability $\nu$, we denote by $T_t(\nu)$ 
the law of $\eta_t$ conditioned 
on $\{\tau>t\}$, with respect to $P_{\nu}$. Thus, for $\v$ continuous and
bounded, $\int \v d T_t(\nu):=E_{\nu}[\v(\eta_t)|\tau>t]$.

Now, from a statistical physics point of view, a relevant issue
is the existence of a limit for $T_t(\nur)$, 
the so-called Yaglom limit, say $\mur$. 
The existence of a Yaglom limit is established by Kesten \cite{k} 
for an irreducible positive recurrent
random walk on $\NN$ with bounded jump size and with $\A=\{0\}$.
It is also established in \cite{af} for the symmetric simple
exclusion process in dimension $d\ge 5$, using strongly the symmetry and
establishing uniform $L^2(\nur)$ bounds 
for $\{dT_t(\nur)/d\nur, t\ge 0\}$. We refer to the introduction of \cite{fk},
for a review of countable Markov chains for which the Yaglom
limit is established. This notion was introduced first by Yaglom 
in 1947 for subcritical branching processes \cite{ya}.

We note that the existence of $\mur$ implies trivialy that there is 
$\lambda(\rho)\in [0,\infty]$
such that for any $s>0$,
\be{eq1.3}
P_{\mur}(\tau>s)=
\lim_{t\to\infty} \frac{P_{\nur}(\tau>t+s)}{P_{\nur}(\tau>t)}
=\exp(-\lambda(\rho) s),
\end{equation} 
and $\lambda(\rho)$ is given by
\be{eq1.3bis}
\lambda(\rho)=-\lim_{t\to\infty}\frac{1}{t} \log\left(P_{\nur}(\tau>t)\right).
\end{equation} 
Thus, right at the outset, one faces three issues.
\begin{itemize} 
\item[(i)] When does
the ratio (\ref{eq1.3}) have a limit? This is linked with a wide area of
investigations
(see e.g. \cite{k,cms,fkm}).
\item[(ii)] Is there a formula for $\l(\rho)$? 
One recognizes in $\l(\rho)$ the logarithm of the
spectral radius of $\LL:L^{\infty}(\nur)
\to L^1(\nur)$ with Dirichlet conditions on $\A$. When $\LL$ is 
a second order elliptic operator on a bounded domain, and when we work
with the sup-norm topology, Donsker and Varadhan~\cite{dv} give
a variational formula for (\ref{eq1.3bis}). 
\item[(iii)] When is $\l(\rho)$
a positive real? 
In other words, what is the right scaling for large deviations for
the occupation time of $\A$. For symmetric simple exclusion, it
is shown in \cite{ar,ad} that $\lambda(\rho)>0$ if and only if $d\ge 3$.
\end{itemize} 

Since $\{T_t,t\ge 0\}$ is a semi-group, the Yaglom limit, when it exists,
is a fixed point of $T_t$ for any $t$. Thus, a preliminary step is
to characterize possible fixed points of $\{T_t\}$,
which are called quasi-stationary measures. We note that
in our context, the Dirac measure on the empty configuration is trivially
a quasi-stationary measure with $\lambda=0$. Thus, by non-trivial
quasi-stationary measure, we mean one corresponding to $\lambda>0$.
Finally, we note that in dynamical systems,
quasi-stationary measures are well studied and
named after Pianigiani and Yorke~\cite{py},
who prove their existence for expanding $C^2$-maps.

Assume that $\mu$ is a probability measure with support in $\A^c$
such that for any $t\ge 0$, $T_t(\mu)=\mu$. By
differentiating this equality at $t=0$, we obtain for
$\v$ in the domain of $\LL$ with $\v|_{\A}=0$
\be{eq1.4}
\int \LL(\v) d\mu=\int \LL(1_{\A^c})d\mu\int \v d\mu.
\end{equation} 
Moreover, assume that $\mu$ is absolutely continuous with respect
to a measure $\nu$, and that $f:=d\mu/d\nu\in L^2(\nu)$. If
$\LL^*$ denotes the adjoint operator in $L^2(\nu)$, 
then $f\in D(\LL^*)$ and $f$ is a nonnegative solution of
\[
1_{\A^c}\LL^*f+\lambda f=0\quad{\rm and}\quad \lambda=
\int -\LL(1_{\A^c})d\mu.
\]
Thus, the problem of quasi-stationary measure for attractive particle
systems is a generalization
of the problem of finding nonnegative eigenvectors, which gave
rise, among others results, to Perron-Frobenius and
Birkhoff-Hopf theorems.
However, such general results cannot be used  in our context, since
neither is the space compact
nor the operator, and since we lack irreducibility conditions.

Equation (\ref{eq1.4}) is the starting point of Ferrari, Kesten,
Mart\'{\i}nez and Picco~\cite{fk}, whose work we describe in 
some details since ours builds upon it.
These authors consider an irreducible,
positive recurrent random walk, $\{X_t, t\ge 0\}$ on $\NN$,
with rates of jump $\{q(i,j),\ i,j\in \NN\}$ 
and study the first time the origin is 
occupied, say $\tau$, when there is $\lambda>0$ and $i\in \NN\bs\{0\}$
such that $E_{i}[\exp(\lambda \tau)]<\infty$.
Assuming that $\mu$ satisfies (\ref{eq1.4}), one obtains for any
$\v$ with $\v(0)=0$
\be{eq1.6}
\sum_{j\not= 0}\sum_{k\not= 0}\left( q(j,k)+q(j,0)\mu(k)\right)
\left( \v(k)-\v(j)\right) \mu(j)=0.
\end{equation} 
Thus, $\mu$ can be thought of as the invariant measure of a 
new random walk, say 
$\{X_t^{\mu}, t\ge 0\}$ on $\NN\bs\{0\}$ with rates
$\{q(j,k)+q(j,0)\mu(k),\ j,k\in \NN\bs\{0\}\}$. 
When $\mu$ is such that $E_{\mu}[\tau]<\infty$,
$X_t^{\mu}$ is positive recurrent and has a unique invariant measure
$\nu$, and this procedure defines a map $\mu\mapsto \Phi(\mu)=\nu$. Thus, 
the problem reduces to finding fixed points of $\Phi$. They notice
also that $X_t^{\mu}$ can be built from the walk $X_t$, by starting
it afresh from a random site drawn from $\mu$, each time $X_t$ hits 0.
Then, using this renewal representation,  
an expression of $\Phi(\mu)$ is obtained (see equation
(2.4) of \cite{fk})
\be{eq1.7}
\Phi(\mu)=\frac{1}{E_{\mu}[\tau]}\int_0^{\infty}\!\!T_t(\mu)P_{\mu}(\tau>t)dt.
\end{equation} 
In our case, equation (\ref{eq1.4}) cannot be interpreted in terms of
$\mu$ being the stationary measure
of a familiar process. Nevertheless, the Laplace-like transform
(\ref{eq1.7}) is a well defined map.
It was observed in \cite{cmm} that as soon as
$E_{\mu}[\tau]<\infty$,
$\mu$ is quasi-stationary if and only if $\Phi(\mu)=\mu$.

In \cite{fk}, the authors study the sequence of
iterates $\{\Phi^{n}(\delta_i)\}$
for $i\in \NN\bs\{0\}$. They show that this sequence is tight, 
and that any limit point belongs to $\M_{\lambda}$,
the subspace of probability measures under which $\tau$ is an
exponential time of parameter
\[
\lambda=-\lim_{t\to\infty} \frac{1}{t}\log\left(P_{\delta_i}(\tau>t)\right)>0.
\]
Then, the facts that $\Phi(\M_{\lambda})\subset \M_{\lambda}$ 
and $\Phi$ is continuous
on the compact set $\M_{\lambda}$,  imply that $\Phi$ has a fixed
point in $\M_{\lambda}$.

Though the irreducibility assumption no longer holds for attractive
particle systems on $\ZZ^d$, we show
that $\{\Phi^{n}(\nur)\}$ is tight through the a priori bounds
$\Phi^{n}(\nur)\prec \nur$, where $\prec$ denotes stochastic domination.
These bounds permit to prove that as soon as $\lambda(\rho)>0$, 
$\tau$ is an exponential time of parameter $\lambda(\rho)>0$,
under any limit point of the iterates sequence. 
We establish that $\lambda(\rho)>0$ in any dimensions for 
zero range processes,
whereas $\lambda(\rho)>0$ is only proved to hold in dimensions larger
or equal than 3 for exclusion processes.

Once $\lambda(\rho)>0$ holds, we show that any
limit point of the Cesaro mean $(\Phi(\nur)+\dots+\Phi^n(\nur))/n$
is quasi-stationary. It is useful to have a sequence converging
to a quasi-stationary measure. Indeed, through a priori bounds,
one gets regularity of the limiting quasi-stationary measure.
For instance, for zero range processes, we can show that 
in dimensions $d\ge 3$, quasi-stationary
measures obtained as Cesaro limits have a density with respect to $\nur$
which is in any $L^p(\nur)$ for $p\ge 1$. In this way, we establish 
the existence of a Dirichlet eigenvector, say $f\in D(\LL^*)$ with
\[
\forall \eta\not\in \A,\quad \LL^* f(\eta)+\lambda(\rho) f(\eta)=0,
\qquad{\rm and}\qquad f|_{\A}=0.
\]
This in turn gives estimates for $P_{\nur}(\tau>t)$ improving
on (\ref{eq1.3bis}).

Finally, we remark that
it could have seemed that a natural way to prove existence of quasi-stationary
measures for our particle systems on $\ZZ^d$, would have been to work first
with finite dimensions approximations, where we can rely on Perron-Frobenius
theory. This strategy fails as is shown on a simple example 
in section~\ref{example}.

\section{Notations and Results.}
\label{notations}
We consider $\NN^{\ZZ^d}$ with the product topology. The local
events are the elements of the union of all $\sigma$-algebras
$\sigma\{\eta(i),i\in \L\}$ over $\L$ finite subset of $\ZZ^d$.
We start by recalling the definition of the
``processus des misanthropes"\cite{co}.
The rates $\{p(i,j),\ i,j\in \ZZ^d\}$ satisfy
\ba{eq1.1bis}
(i)&& p(i,j) \geq 0, \quad \sum_{i\in \ZZ^d} p(0,i)=1. \cr
(ii)&&p(i,j)=p(0,j-i)\quad\text{(translation invariance)}.\cr
(iii)&& p(i,j)=0\text{ if }|i-j|>R\text{ for some fixed }R
\quad\text{(finite range)}.\cr
(iv)&& {\rm If\ }p_s(i,j)=p(i,j)+p(j,i),\text{ then}\quad
\forall i\in \ZZ^d,\ \exists n,\quad p_s^{(n)}(0,i)>0
\quad\text{(irreducibility)}.\cr
(v)&& \sum_{i\in \ZZ^d} ip(0,i)\not= 0 \quad\text{(drift)}.
\ea
Let $b:\NN\times \NN\to [0,\infty)$ be a function with 
\ba{def.b}
(i)&& b(0,.)\equiv 0\cr
(ii)&& n\mapsto b(n,m) \quad\text{ is nondecreasing for each }m\cr
(iii)&& m\mapsto b(n,m) \quad\text{ is nonincreasing for each }n\cr
(iv)&& b(n,m)-b(m,n)=b(n,0)-b(m,0),\quad \forall n,m\ge 1\cr
(v)&& \D:=\sup_{n}\left(b(n+1,0)-b(n,0)\right)<\infty.
\ea
As in~\cite{andjel}, a Feller process can be constructed on
\[
\Omega=\{\eta:\ \sum_{i\in \ZZ^d} e^{-a|i|}\eta(i)<\infty,
\text{ for some }a>0\},
\]
with generator acting on a core of local functions as
\be{def.g}
\LL \v(\eta):=\sum_{i,j\in \ZZ^d}p(i,j)b(\eta(i),\eta(j))\left(
\v(\eta^i_j)-\v(\eta)\right),
\end{equation} 
where $\eta^i_j(k)=\eta(k)$ if $k\not\in \{i,j\}$, 
$\eta^i_j(i)=\eta(i)-1$, and $\eta^i_j(j)=\eta(j)+1$.

Let $g:\NN\to [0,\infty)$ satisfy (\ref{eq1.0}), and $g(1)=1$.
For any $\g\in [0,\sup_k g(k)[$, we define
a probability $\theta_{\g}$  on $\NN$, by
\be{def.m}
\theta_{\g}(0)=1/Z(\g),\quad\text{ and when }n\not= 0,\quad
\theta_{\g}(n)=\frac{1}{Z(\g)}\frac{\g^n}{g(1)\dots g(n)},
\end{equation} 
where $Z(\g)$ is the normalizing factor. If we set
$\Upsilon(\g)=\sum_{n=1}^{\infty} n\theta_{\g}(n)$, then
$\Upsilon:[0,\sup_k g(k)[\to [0,\infty[$ is increasing.
Let $\gamma: [0,\sup_{\gamma}\Upsilon(\gamma) ) \to [0,\sup_k  g(k) )$ 
be the inverse of $\Upsilon$, and let $\nur$ be the product
probability with marginal law $\theta_{\g(\rho)}$. Thus, we have
\be{def.nur}
\forall i\in \ZZ^d,\quad\int \eta(i)d\nur=\rho,\quad{\rm and}\quad
\int g(\eta(i))d\nur=\g(\rho).
\end{equation}
For a function $b$ satisfying (\ref{def.b}), 
we assume there is  $g$ as above, with
$ b(n,m-1)g(m)=b(m,n-1)g(n)$, which together with (\ref{def.b} (iv))
and (\ref{eq1.1bis} (i)), 
imply that $\{\nur, \rho\in [0, \sup_{\gamma}\Upsilon(\gamma) )\}$ 
are invariant with respect to $\LL$.

Now, if we choose $b(n,m)=g(n)$, we obtain the zero range process.
We describe a way of realizing this process, in case like ours,
where the labelling of particles is innocuous. We start with 
an initial configuration $\eta\in \Omega$. We label arbitrarily
particles on each site $i$ from 1 to $\eta(i)$. We associate
to each particle a path $\{S_n,\ n\in \NN\}$, paths being drawn
independently from those of a random walk with rates $\{p(i,j)\}$.
Then, a particle labelled $k$ at site $i$ jumps with rate
$g(k)-g(k-1)$. If it jumps on site $j$ it gets the last label. Also,
the remaining particles at site $i$ are relabelled from 1 to $\eta(i)-1$.
Now, as $\D:=\sup_{k>1}\left(g(k)-g(k-1)\right)<\infty$, we can dominate
the Poisson clocks with independent Poisson clocks of intensity $\D$, so that 
each particle is coupled with a random walk wandering faster on the same path. 

If we restrict the process to $\{0,1\}^{\ZZ^d}$, and choose
$b(n,m)=1$ if $n=1,m=0$ and $b(n,m)=0$ otherwise, we obtain
the exclusion process. The measure $\nur$ is then a product Bernoulli measure.

The semi-group $\{S_t\}$ generated by $\LL$ extends to a Markov semi-group
on $L^2(\nur)$, and its generator is the closure of $\LL$ to 
$L^2(\nur)$ (see the proof of Prop. 4.1 of~\cite{liggett}).
We can consider also the adjoint (or time-reversed)
of $\LL$ in $L^2(\nur)$, as acting on local functions $\v$ and $\psi$ by
\be{def.adjoint}
\int \LL^*(\v) \psi d\nur:=\int \v\LL(\psi)d\nur.
\end{equation}
With our hypothesis, $\LL^*$ is again the generator of a 
``processus des misanthropes''  on $\Omega$, 
with the same functions $b$ and $g$,
but with $p^*(i,j):=p(j,i)$ (see e.g. \cite{ba}). We denote
by $\{S_t^*\}$ the associated semi-group, and by $P^*_{\eta}$ the 
associated Feller process with initial configuration $\eta\in \Omega$.

For convenience, we fix an integer $k$ and $\L$ a finite subset of $\ZZ^d$,
and set $\A:=\{\eta:\ \sum_{i \in \L}\eta(i)>k\}$. Needless to emphasize that
we will always consider a density $\rho$ such that $\nur(\A^c)>0$.
We denote by $\bar \LL:=1_{\A^c}\LL$ and
$\{\bar S_t,t\ge 0\}$, respectively the generator and associated semi-group
for the process killed on $\A$. A core of $\bar \LL$ consists of local
functions vanishing on $\A$.

For $\eta,\xi \in \Omega$, we say that $\eta \leq \xi$
if $\eta(i) \leq \xi(i)$ for all $i \in \ZZ^d$. Monotonicity of
functions from $\Omega$ to $\RR$ is meant with this partial
order; in particular, we say that $A \subset \Omega$ is
increasing if $1_A$ is increasing. Finally, for
given probability measures $\nu,\mu$ on $\Omega$, we say that
$\nu \prec \mu$ if $\int f d\nu \leq \int f d\mu$ for every
increasing function $f$. We recall that the ``processus des misanthropes''
is an attractive process, i.e. 
there is a coupling such that $P_{\eta,\zeta}(\eta_t\le \zeta_t,
\forall t)=1$ whenever $\eta\le \zeta$.

Since $\A$ is an increasing local event, attractiveness
implies that for any $t\ge 0$,
both $P_{\eta}(\tau>t )$ and $P^*_{\eta}(\tau>t )$ are decreasing in 
$\eta$. As our product measure satisfies FKG's inequality, we have
\be{eq.sub}
P_{\nur}(\tau>t+s)=\int \bar S_{t+s}(1_{\A^c})d\nur=
\int \bar S_{t}(1_{\A^c})\bar S_{s}^*(1_{\A^c})d\nur\ge
P_{\nur}(\tau>t)P_{\nur}(\tau>s).
\end{equation}
Also it is easy to see that $\nur(\A^c)>0$ implies that for
any $t\ge 0$, $P_{\nur}(\tau>t)>0$
(this is true for short time by continuity, and one then uses (\ref{eq.sub})
to extend it to any time).
Thus, (\ref{eq.sub}) and $P_{\nur}(\tau>t)>0$ justify 
the existence of the limit $\l(\rho)<\infty$ in (\ref{eq1.3bis}).

A key, though elementary, observation of \cite{fk,cmm} is as follows.
\bl{lem1} Let $\mu$ be such that $E_{\mu}[\tau]<\infty$. Then, 
$\mu$ is quasi-stationary if and only if $\Phi(\mu)=\mu$.
\el
Indeed, if $\mu$ is quasi-stationary, 
then it is obvious that $\Phi(\mu)=\mu$. Conversely, for any $\v\in \C_b$
\[
\int \bar S_s(\v)d\mu=\frac{1}{E_{\mu}[\tau]}\int_0^{\infty}
\int\bar S_t(\bar S_s(\v))d\mu dt=
\frac{1}{E_{\mu}[\tau]}\int_s^{\infty} \int\bar S_t(\v)d\mu dt,
\]
which implies that 
\[
\int \bar S_s(\v)d\mu=\exp(-\frac{s}{E_{\mu}[\tau]})\int \v d\mu.
\]
Now, a key a priori bound relies on the notion of stochastic domination.
\bl{lem2} If $\Phi^n$ denotes the $n$-th iterate of $\Phi$, 
then $\Phi^n(\nur)\prec \nur$. Also, $\{\Phi^n(\nur)\}$ is tight.
\el
This allows us to prove a result analogous to Lemma 3.2 of \cite{fk}.
\bl{lem3} If $\lambda(\rho)\in ]0,\infty[$, then for any integer $k\ge 1$
\[
\lim_{n\to\infty} \int \tau^kd\Phi^n(\nur)=\frac{k!}{\lambda(\rho)^k}.
\]
Moreover, for any $s\ge 0$
\be{eq1.8}
\lim_{n\to\infty}P_{\Phi^n(\nur)}(\tau>s)=
\exp(-\lambda(\rho) s).
\end{equation}
\el
If we set $\bar \nu_n:=\frac{1}{n}(\Phi(\nur)+\dots+\Phi^n(\nur))$,
then our existence result reads.
\bt{the1} Assume that $\l(\rho)>0$. Then, any limit point along a 
subsequence of $\{ \bar \nu_n,\ n\in \NN\}$
is a quasi-stationary measure corresponding to $\lambda(\rho)$. 
\et
We prove Lemmas~\ref{lem2}, \ref{lem3} and Theorem~\ref{the1} in
section~\ref{existence}.
We give now conditions under which $\l(\rho)>0$.
Note that in the symmetric case, \cite{ad} established the following
stronger result using spectral representation.
\be{eq1.9}
\lim_{u\to\infty} \frac{P_{\nur}(\tau>u+s)}{P_{\nur}(\tau>u)}=
e^{-\lambda_{s}(\rho) s}\quad{\rm with}\quad 
\lambda_{s}(\rho)= \inf\Bigl\{
\frac{-\int f \LL fd\nur}{\int f^2d\nur}:\ f\in D(\LL),\ f|_{\A}=0\Bigr\}.
\end{equation}
It was established in \cite{ad} that 
for the symmetric exclusion process $\lambda_{s}(\rho)>0$ for $d\ge 3$, and
that $\lambda_{s}(\rho)=0$ for $d=1$ and $d=2$.
Using the classical bound $\lambda(\rho)\ge \lambda_{s}(\rho)$ 
(see e.g. \cite{reza} Lemma 4.1), we have 
\bl{lem.known} For the exclusion process in $d\ge 3$, $\lambda(\rho)$
given by (\ref{eq1.3bis}) is positive.
\el
For zero range processes, 
we prove in section~\ref{regularity} the following results.
\bl{lem4} For zero range processes in any dimensions, $\lambda(\rho)>0$.
\el
Moreover, we have the following regularity result.
\bp{lem5} For zero range processes in $d\ge 3$, 
any limit points along a subsequence of $\{\bar \nu_n\}$,
say $\mur$, is absolutely continuous with respect
to $\nur$ and $f:=d\mur/d\nur \in L^p(\nur)$ for any $p\ge 1$.
Thus, $f$ is in the domain of $\bar \LL^*$ and
\be{eq.eigen}
\bar\LL^*f+\lambda(\rho) f=0,\quad{\rm a.s.}-\nur.
\end{equation}
\ep
As a consequence of the existence of an eigenvector of (\ref{eq.eigen})
in $L^p(\nur)$ for $p\ge 1$, we have estimates for the hitting time.
\bc{cor1} For zero range processes in $d\ge 3$,
let $f$ be a solution of (\ref{eq.eigen}) and $g$ a solution of the adjoint
eigenvector equation. Then, $\int fgd\nur$ is finite
and positive, and for any time $t$
\be{eq.hitting}
\exp(-H(\tilde \nur,\nur))\le \frac{P_{\nur}(\tau>t)}{\exp(-\lambda(\rho) t)}
\le 1,
\end{equation}
with
\[
d\tilde \nur=\frac{fgd\nur}{\int fgd\nur},\quad
\text{and}\quad H(\tilde\nur,\nur)
=\int \log(\frac{d\tilde \nur}{d\nur})d\tilde \nur<\infty.
\]
\ec

Finally, in section~\ref{example} we see,
on the totally asymmetric simple exclusion process, why the finite dimensional
approximation of our problem yields `wrong' results.

\section{Existence.}
\label{existence}
We begin with some useful expressions for the iterates $\nu_n:=\Phi^n(\nur)$.
If $\l(\rho)>0$, then $\forall n \in \NN$, 
$\int_0^{\infty}u^n P_{\nur}(\tau>u)du$ is finite,
and it follows easily by induction that
\be{eq.nun}
\int \v d\nu_n=\frac{\idotsint_0^{\infty}\int\bar S_{t_1+\dots+t_n}(\v)d
\nur\prod_{i=1}^ndt_i}{\idotsint_0^{\infty}\int\bar S_{t_1+\dots+t_n}
(1_{\A^c})\nur\prod_{i=1}^ndt_i}=
\frac{\int_0^{\infty}u^{n-1}\int \bar S_u(\v) d\nur du}{
\int_0^{\infty}u^{n-1}\int \bar S_u(1_{\A^c})d\nur du}.
\end{equation}
Applying this expression to $\v=\bar S_t(1_{\A^c})$ yields
\[
P_{\nu_n}(\tau>t)=
\frac{\int_0^{\infty}u^{n-1}P_{\nur}(\tau>t+u)du}{
\int_0^{\infty}u^{n-1}P_{\nur}(\tau>u)du}.
\]
Integrating over $t$, we obtain
\be{eq.temps}
E_{\nu_n}[\tau]=\frac{1}{n}
\frac{\int_0^{\infty}u^{n}P_{\nur}(\tau>u)du}{
\int_0^{\infty}u^{n-1}P_{\nur}(\tau>u)du}
=\frac{E_{\nur}[\tau^{n+1}]}{(n+1)E_{\nur}[\tau^{n}]}.
\end{equation}
\noindent
{\bf Proof of Lemma~\ref{lem2}.}
Let $\v$ be a nondecreasing bounded function, then
\[
\int \bar S_u \v \, d\nur=\int E_{\eta}[\v(\eta_u)1_{\{\tau>u\}}] \, d\nur=
\int \v(\eta)\bar S_u^*(1_{\A^c})(\eta) \, d\nur
\]
Now, we note that $\eta\mapsto\bar S_u^*1_{\A^c}(\eta)$ is nonincreasing. 
By FKG's inequality, we thus have 
\[
\int \bar S_u \v \, d\nur\le \int \v \, d\nur \int \bar S_u(1_{\A^c}) \,
d\nur \,.
\]
This implies by (\ref{eq.nun}) that 
$
\int \v \, d\nu_n \le
\int \v \, d\nur$.
Consider now compact subsets of $\NN^{\ZZ^d}$ of the type 
$K_{(k_i)}=\{\eta:\, \forall i \in \ZZ^d,
\eta_i \leq k_i\}$. Since these compacts are decreasing, we have
$\inf_n \nu_n(K_{(k_i)}) \geq \nur(K_{(k_i)})$. Moreover, 
for all $\epsilon >0$,
a good choice of the sequence $(k_i)$ ensures that   
$\nur(K_{(k_i)}) \geq 1-\epsilon$, and tightness follows.
\qed

\noindent
{\bf Proof of Lemma~\ref{lem3}.}
The argument follows closely~\cite{fk} (proofs of Lemma 3.2,
Proposition 3.3 and Theorem 4.1), the main difference being that
we replace irreducibility by stochastic domination.
If $\nu_n=\Phi^n(\nur)$, then we show in three steps
that $ \lim E_{\nu_n}[\tau]=1/\lambda(\rho)$.

\noindent
\underline{Step 1}: We first prove that 
\be{eq1.17}
\underline\lim E_{\nu_n}[\tau]=1/\lambda(\rho)\quad{\rm and}\quad
P_{\nur}(\tau>t)\le \exp(-\lambda(\rho)t).
\end{equation}
As in Proposition 3.3 of \cite{fk}, if
\[
\frac{1}{\lambda_{\infty}}=\underline{\lim} E_{\nu_n}[\tau],
\quad{\rm then}\quad \lambda_{\infty}\ge \lambda(\rho),
\]
and there is a subsequence $\{n_k\}$ such that
\[
\forall t>0,\qquad \lim_{k\to\infty}P_{\nu_{n_k}}(\tau>t)=
\exp(-\lambda_{\infty}t).
\]
The inequality $\lambda_{\infty}\le \lambda(\rho)$ follows 
after observing that as $\eta\mapsto P_{\eta}(\tau>t)$ is decreasing,
and as $\nu_n\prec \nur$, we have
$ P_{\nu_{n_k}}(\tau>t)\ge P_{\nur}(\tau>t)$. Thus,
\be{eq1.16}
\exp(-\lambda_{\infty}t)=\lim_{k\to\infty}
P_{\nu_{n_k}}(\tau>t)\ge P_{\nur}(\tau>t).
\end{equation}
This establishes that $\lambda_{\infty}=\lambda(\rho)$ and (\ref{eq1.17}).

\noindent
\underline{Step 2}: We show that 
\be{eq1.18}
\lim_{n\to\infty} \left( \frac{E_{\nur}[\tau^n]}{n!}\right)^{1/n}=
\frac{1}{\lambda(\rho)}.
\end{equation}
First, by step 1,
\be{eq1.19}
E_{\nur}[\tau^n]=\int_0^{\infty} nu^{n-1} P_{\nur}(\tau>u)du\le
\int_0^{\infty} nu^{n-1} \exp(-\lambda(\rho)u) du=
\frac{n!}{\lambda(\rho)^n}.
\end{equation}
If we set $v_n=E_{\nur}[\tau^n]/n!$, we have then $\lim\sup v_n^{1/n}\le 1/
\lambda(\rho)$. Now, by (\ref{eq.temps}), $E_{\nu_n}[\tau]=v_{n+1}/v_n$. 
Since $\underline\lim E_{\nu_n}[\tau]=\frac{1}{\lambda(\rho)}$, 
it follows that 
\be{eq1.20}
\forall \e\in ]0,1/\l(\rho)[,\ \exists n_0,\ \forall n\ge n_0, \qquad
v_n\ge v_{n_0}\left(\frac{1}{\lambda(\rho)}-\e\right)^{n-n_0}.
\end{equation}
Thus, for any $\e>0$, $\underline\lim\  v_n^{1/n}\ge 1/\lambda(\rho)-\e$,
and this concludes step 2.

\noindent 
\underline{Step 3}: 
We show that $\overline\lim E_{\nu_n}[\tau]\le 1/\lambda(\rho)$
by following the proof of Theorem 4.1 of \cite{fk}. We omit the argument
here.

Finally, as in \cite{fk}, it is now easy to conclude that
for any integer $k\ge 1$ and $s>0$
\[
E_{\nu_n}[\tau^k]=k!\prod_{j=1}^k E_{\nu_{n+j+1}}[\tau]\to \frac{k!}
{\lambda(\rho)^k},\quad{\rm and}\quad
P_{\nu_n}(\tau>s)\to e^{-\lambda(\rho)s}.
\]
\qed

\noindent
{\bf Proof of Theorem~\ref{the1}.}
For any integer $n$, set $\bar \nu_n=(\Phi(\nur)+\dots+\Phi^n(\nur))/n$.
Note that from Lemma~\ref{lem2} and Lemma~\ref{lem3}, we have that
\be{eq3.1}
\bar \nu_n\prec \nur,\quad
E_{\bar \nu_n}[\tau^k]
\begin{array}[t]{c}
\longrightarrow \\[-2mm] \scriptstyle{ n \rightarrow \infty}
\end{array} 
\frac{k!}{\lambda(\rho)^k},\quad{\rm and}\quad
P_{\bar \nu_n}(\tau>t)
\begin{array}[t]{c}
\longrightarrow \\[-2mm] \scriptstyle{n \rightarrow \infty}
\end{array}
\exp(-\lambda(\rho)t).
\end{equation}
Thus, $\{\bar \nu_n\}$ is tight and let $\mu$ be a limit point along
subsequence $\{\bar\nu_{n_k}\}$.
As $\A^c$ is local and $\bar S_t$ is Feller, (\ref{eq3.1}) implies that
\be{eq3.2}
P_{\mu}(\tau>t)=\lim_{k\to\infty}\!\! P_{\bar \nu_{n_k}}(\tau>t)=
e^{-\l(\rho)t}.
\end{equation}
We now check that $\Phi(\mu)=\mu$, or in other words, that
for $\v$ continuous and bounded
\be{eq3.3}
\lambda(\rho)\int_0^{\infty}\int \bar S_t\v d\mu dt=\int \v d\mu.
\end{equation}
Now, for all $t\ge 0$, the integrable bound
\[
|\int \bar S_t\v d\bar\nu_{n_k}|\le |\v|_{\infty} P_{\bar \nu_{n_k}}(\tau>t)
\le |\v|_{\infty} \left( 1\wedge \frac{\sup_n E_{\bar \nu_n}[\tau^2]}{t^2}
\right),
\]
and $\lim_k\int \bar S_t\v d\bar\nu_{n_k}=\int \bar S_t\v d\mu$ imply,
by dominated convergence, that
\be{eq3.3bis}
\lim_k\int_0^{\infty}\left(\int \bar S_t\v d\bar\nu_{n_k}\right)dt=
\int_0^{\infty}\left(\int \bar S_t\v d\mu\right)dt.
\end{equation}
However, by definition of the iterates
\[
\int \v d\nu_{k+1}=\frac{\int \int_0^{\infty} \bar S_t(\v) dt d\nu_k}
{E_{\nu_k}[\tau]}.
\]
Thus,
\be{eq3.4}
\int \int_0^{\infty}\bar S_t\v dtd\bar \nu_{n_k}=
\frac{1}{n_k}\sum_{i=1}^{n_k} E_{\nu_{i}}[\tau] \int \v d\nu_{i+1}
\longrightarrow \frac{1}{\lambda(\rho)}\int \v d\mu.
\end{equation}
The result follows by (\ref{eq3.3bis}) and (\ref{eq3.4}).
\qed

\section{Positivity of $\lambda(\rho)$ and regularity.}
\label{regularity}
Let $\AA_i:\Omega\to\Omega$ with $\AA_i\eta(k)=\eta(k)+\d_{i,k}$.
For any continuous and bounded function $\v$, we have
\be{eq3.9}
\int g(\eta_i)\v\ d\nur=\g(\rho) \int \AA_i(\v)d\nur.
\end{equation}
Note also that as $k\D\ge g(k)$, we have
\be{eq3.9bis}
\int \eta_i\v\ d\nur\ge\frac{\g(\rho)}{\D} \int \AA_i(\v)d\nur.
\end{equation}

\noindent
{\bf Proof of Lemma~\ref{lem4}.}
We prove that $\Pn\le \exp(-\lambda t)$ for $\lambda>0$, by showing
that
\be{eq3.5}
-\frac{d\Pn}{dt}=-\int \bar S_t(\bar \LL 1_{\A^c})d\nur\ge
\l \int \bar S_t(1_{\A^c})d\nur.
\end{equation}
Now,
\be{eq3.6}
-\bar \LL 1_{\A^c}(\eta)=\sum_{i\not\in \L}\sum_{j\in \L}
p(i,j)g(\eta_i)1_{\{\eta\not\in \A,\eta^i_j\in \A\}}.
\end{equation}
We set $\p \A:=\{\eta:\ \sum_{\L}\eta(i)=k\}$ and note that 
since $g(0)=0$, for any $i\not\in \L$ and any $j\in \L$,
$
g(\eta_i)1_{\p\A}=g(\eta_i)1_{\{\eta\not\in \A,\eta^i_j\in \A\}}$.
Hence,
\[
\begin{array}{ll}
-\int \bar S_t(\bar \LL 1_{\A^c}) \, d\nur
& = - \int \bar \LL1_{\A^c} \,   \Pe \, d \nur
= \sum_{i\not\in \L , j\in \L} p(i,j) \int_{\p \A}  g(\eta_i)\Pe \, d\nur
\\
& =  \g(\rho) \sum_{i\not\in \L , j\in \L} p(i,j)  \int_{\p \A}  \Pi d\nur \, ,
\end{array}
\]
where we have used (\ref{eq3.9}) and the fact that
$\p \A$ is independent of $\eta_i$ for $i\not\in \L$.

Since $\{(i,j) \in \L^c \times \L, \text{s.t. } p(i,j) > 0\}$ is finite, 
we have
now to prove that $\forall i \notin \L$, $\exists \lambda_i > 0$ such that
\[
\int_{\p \A}\Pi d\nur\ge \lambda_i \int \Pe d\nur.
\]
This will be done in three steps.

\noindent 
\underline{Step 1}: 
We show that for $i\not\in \L$, there is $\e_i>0$ such that
\be{eq3.7}
\Pi \ge \e_i \Pe.
\end{equation}
We need to couple two trajectories, say $\{\eta_t,\zeta_t\}$ differing
by a particle at $i$ at time 0, i.e. $\zeta_0=\AA_i \eta_0$.
We describe a basic coupling. We tag the additional particle at $i$,
and call its trajectory $\{X(i,t),t>0\}$. It follows the path 
$\{S_n, n\in \NN\}$ of a random walk with rates $p(.,.)$,
and jumps at the time-marks of an $\eta$-dependent
Poisson clock: at time $t$, its intensity is $g(\eta_t(X(i,t))+1)-
g(\eta_t(X(i,t)))$. With this labelling, the motion of the additional
particle does not perturb the $\eta$-particles. Thus, we call the
additional particle a $2^{nd}$-class particle.
As $\Delta:=\sup(g(k+1)-g(k))<\infty$,
we can couple $\{X(i,t),t>0\}$, with $\{\tilde X(i,t),t>0\}$ which follows
the same path $\{S_n, n\in \NN\}$, but with a Poisson clock of intensity
$\Delta$ which dominates the clock of $\{X(i,t),t>0\}$. Thus,
\be{eq3.8}
S(\L^c)=\inf\{t:\ X(i,t)\in \L\}\ge
\tilde S(\L^c)=\inf\{t:\ \tilde X(i,t)\in \L\}.
\end{equation}
and under our coupling, we have that 
$\{S(\L^c)<\infty\}\subset \{\tilde S(\L^c)<\infty\}\subset 
\{S_n\in \L, n\in \NN\}$. Therefore,
\be{gradient}
\begin{array}{ll}
0\le \Pe-\Pi
&\le P^*_{\eta}(\tau(\eta_.)>t,\ \tau(\zeta_.)\le t) \\
&\le P^*_{\eta}(\tau(\eta_.)>t,\ S(\L^c)<\infty)  \\
& \le P^*_{\eta}(\tau(\eta_.)>t,\tilde S(\L^c)<\infty) \\
&\le\PP_i(S_n\in \L,\ n\in \NN) \Pe \, .
\end{array}
\end{equation} 
Now, as the walk is transient, 
$\e_i:=\PP_i(S_n\not\in \L ,\forall n\in \NN)>0$,
so that (\ref{eq3.7}) holds. 

\noindent 
\underline{Step 2}: 
It remains now to show that $\int_{\p \A} \Pe \, d\nur \ge \lambda \int
\Pe \, d\nur$ for some $\lambda >0$. This would be easily done by FKG
inequality, if $\p \A$ was a decreasing event, which is not the case.
However, $\A_0 :=\{\eta: \, \sum_{i\in \L} \eta(i)=0\}$ is a decreasing event,
 and the idea is to compare $\int_{\p \A} \Pe \, d\nur$ with 
$\int_{\A_0} \Pe \, d\nur$. 
To this end, we are going to compare  $P^*_{\eta}(\tau > t)$ for
$\eta \in \p \A$, with $P^*_{\AA_j^{-1} \eta}(\tau > t)$ for $j \in \L$,
so that  we consider now the case where the $2^{nd}$-class particle
is initially in $j\in \L$. We will ensure that,
uniformly in $\eta\in \p \A$, there is a positive
probability that the $2^{nd}$-class particle escapes $\L$ within a small
time $\d>0$. If the $2^{nd}$-class particle
finds itself on a site with $k$ particles, it jumps with
rate $\Delta_k:=g(k+1)-g(k)$. We have $\Delta_1>0$, but could
very well have $\Delta_k=0$ for $k>1$. Thus, 
the $2^{nd}$-class particle can move for sure only when on an empty site.
As in Step 1, we have a coupling $(\eta_.,\zeta_.)$, where $\zeta_0=
\AA_j\eta_0$. For convenience, we use the notation
$P_{\eta,j}$ instead of $P_{\zeta}$.

Thus, we impose on the $\eta$-particles starting on $\L$ the following
constraints: 
\begin{itemize} 
\item[(i)] they do not escape from $\L$ during $[0,\d]$;
\item[(ii)] they empty one `path' joining $j$ with $\p \L$ during $[0,\d/3]$
while freezing the $2^{nd}$-class particle;
\item[(iii)] we freeze their motion during $[\d/3,2\d/3]$ while forcing
the $2^{nd}$-class particle to escape $\L$;
\item[(iv)] we force the $\eta$-
particles to go back to their initial configuration during $]2\d/3,\d]$.
\end{itemize} 
More precisely, we let $\Gamma:=\{j_1,\dots,j_n\}$ be a shortest path
linking $j$ to $\L^c$, that is
\[
j_1=j,\ j_2,\dots,j_{n-1}\in \L,\text{ and }j_n\not\in \L,
\quad\text{and\ } p(j_k,j_{k+1})>0 \text{ for } k<n.
\]
We note $i_j:=j_n$ the extremity of $\G$, 
and for a subset $A$ of $\ZZ^d$, we call $\s(A)$ the first time
that an $\eta$-particle initially in $A$ exits $A$. Also, let
\[
D_{\L}:=\{\eta:\ \eta(j_k)=0 \text{ for } k=1,\dots,n-1\}\cap \p \A.
\]
Now, we say that $(\eta_.,X(j,.))\in \F_{j,i_j} [0,\d]$ if 
\begin{itemize}
\item[(i)] $\s(\L)(\eta_.)>\d$;
\item[(ii)] on $[0,\d/3]\ X(j,.)=j$ and $\eta_{\d/3}\in D_{\L}$;
\item[(iii)] on $[\d/3,2\d/3]$, $\eta_.|_{\L}=\eta_{\d/3}|_{\L}$, 
$X(j,.)$ reaches $i_j$ before $2\d/3$ along $\G$, and stays still;
\item[(iv)] on $[2\d/3,\d]\ X(j,.)=i_j$,  and $\eta_.|_{\L}
=\eta_{\d-t}|_{\L}$.
\end{itemize} 
We call $\tilde \F_{i_j,j} [0,\d]$ the time reversed event
\[
\{(\eta_.,X(i,.))\in \tilde \F_{i_j,j} [0,\d]\}:=
\{(\eta_{\d-.},X(j,\d-.))\in \F_{j,i_j} [0,\d]\}.
\]
It is plain that 
\be{eq3.11}
\lambda_1:=\inf_{\eta: \sum_{i \in \L}\eta(i)\le k} \, \, \inf_{j\in \L}
P_{\eta,j}^*(\F_{j,i_j}[0,\d])>0.
\end{equation}
We prove in this step that there is $\lambda_2>0$ such that for $\eta$ such
that $\sum_{i \in \L} \eta(i) \leq k-1$, 
\be{eq3.12}
\Pj=P^*_{\eta,j}(\tau(\zeta_.)>t)
\ge \lambda_2 P^*_{\eta,j}(\tau(\eta_.)>t, \s(\L^c)>\d, \F_{j,i_j}[0,\d]).
\end{equation}
>From the instant $\d$, we couple through our basic coupling,
the $2^{nd}$-class particle with a random walk whose
Poisson clock has intensity $\Delta$, so that
\be{eq3.13}
\{\tilde S(\L^c)\circ \theta_{\d} =\infty\}\subset 
\{S(\L^c)\circ \theta_{\d}=\infty\}.
\end{equation}
Note that if  particles from outside $\L$, do not enter $\L$ during time
$[0,\d]$,  if the $2^{nd}$-class particle exits $\L$ before $\d$, not to ever
enter again, and if $\{\tau(\eta_.)>t\}$, then $\{\tau(\zeta_.)>t\}$.
In other words,
\be{eq3.14}
\{\tau(\eta_.)>t\}\cap \{\s(\L^c)>\d\}\cap \F_{j,i_j}[0,\d]
\cap \{S(\L^c)\circ \theta_{\d}=\infty\}
\quad\subset\quad \{\tau(\zeta_.)>t\}.
\end{equation}
Thus, by conditioning on $\sigma\{\zeta_s,s\le \d\}$
\begin{align}
P^*_{\eta,j}(\tau(\zeta_.)>t)\ge 
&P^*_{\eta,j}(\tau(\eta_.)>t,\s(\L^c)>\d,\F_{j,i_j}[0,\d]
,S(\L^c)\circ \theta_{\d}=\infty) \notag\\
\ge & P^*_{\eta,j}(\tau\circ \theta_{\d}(\eta_.)>t,\s(\L^c)>\d,
	\F_{j,i_j}[0,\d],
\tilde S(\L^c)\circ \theta_{\d} =\infty) \notag\\
\ge & E_{\eta,j}^*[1_{\{\sigma(\L^c)>\d,\F_{j,i_j}[0,\d]\}}
P^*_{\eta_{\d},i_j}(\tau(\eta_.)>t-\delta ,\tilde S(\L^c)=\infty)]\notag\\
\ge &\PP_{i_j}(S_n\not\in \L, \forall n\in \NN)
P^*_{\eta,j}(\tau(\eta_.)>t,\s(\L^c)>\d,\F_{j,i_j}[0,\d]).\notag
\end{align}
This is (\ref{eq3.12}), once we recall that $\{S_n\}$ is transient, and
that $\{i_j; \,  j \in \L\}$ is finite.

\noindent 
\underline{Step 3}: We prove the result inductively. We fix one configuration
in $\p \A$: let $\{k_j, j\in \L\}$, be integers such that
\be{eq3.16}
\sum_{j\in \L} k_j=k,\quad{\rm and}\quad
\B:=\{\eta: \eta_j=k_j, j\in \L\}.
\end{equation}
Let $j$ be such that $k_{j}>0$. Then, using (\ref{eq3.9bis})
\begin{align}
\int_{\B}\!\! \Pe d\nur=& 
\int_{\B} \frac{\eta_{j}}{k_{j}} \Pe d\nur(\eta)\notag\\
\ge&\frac{\g(\rho)}{\D k_{j}}\int_{\AA_{j}^{-1}\B}\!\!
\Pj d\nur(\eta)\notag\\
\ge & \frac{\lambda_2\g(\rho)}{\D k_{j}}\int_{\AA_{j}^{-1}\B}\!\!
P^*_{\eta,j}(\tau(\eta_.)>t,\s(\L^c)>\d,\F_{j,i_j}[0,\d])d\nur.\notag
\end{align}
Using the stationarity of $\nur$, and reversing
time on the interval $[0,\d]$, the last integral becomes
\[
\int P_{\eta,i_j}(\tilde \F_{i_j,j}[0,\d],\eta_{\d}\in 
\AA_{j}^{-1}\B, \s(\L^c)>\d) P^*_{\eta}(\tau>t-\d)d\nur(\eta).
\]
Note that in $\{\tilde \F_{i_j,j}[0,\d],\eta_{\d}\in \AA_{j}^{-1}\B, \s(\L^c)>\d
\}$, the particles from inside and outside $\L$ do not interact, 
and that $\tilde \F_{i,j}[0,\d]$ imposes the same initial and final 
configuration for the $\eta$-particles in $\L$, so that
\[
P_{\eta,i_j}(\tilde \F_{i_j,j}[0,\d],\eta_{\d}\in \AA_{j}^{-1}\B, \s(\L^c)>\d)=
1_{\B}(\AA_{j}(\eta))P_{\eta,j}^*(\F_{j,i_j}[0,\d]) P_{\eta}(\s(\L^c)>\d).
\]
Thus, from (\ref{eq3.11}), there is $\tilde \e>0$ such that
\be{eq3.18}
\int_{\B}\!\! \Pe d\nur
\ge\tilde \e \int_{\AA_{j}^{-1}\B}\!\!
P_{\eta}(\s(\L^c)>\d)P^*_{\eta}(\tau>t-\d)d\nur(\eta).
\end{equation}
We iterate the same procedure $k$ times, and end up with $\e>0$ such that
\be{eq3.19}
\int_{\B}\!\! \Pe d\nur\ge \e 
\int_{\prod_{j\in \L}\AA_{j}^{-k_j}\B}\!\!
P_{\eta}(\s(\L^c)>k\d)P^*_{\eta}(\tau>t-k\d)d\nur(\eta).
\end{equation}
Finally, we note that 
\[
\eta\mapsto 1_{\prod_{j\in \L}\AA_{j}^{-k_j}\B}
= 1_{\{\eta:\,  \eta(j)=0, j\in \L\}},
\quad
\eta\mapsto P_{\eta}(\s(\L^c)>k\d),\quad{\rm and}\quad
\eta\mapsto P^*_{\eta}(\tau>t-k\d),
\]
are decreasing functions. Thus, by FKG's inequality
\be{eq3.20}
\int_{\B}\!\! \Pe d\nur\ge \e\nur(\{\eta: \eta(j)=0, j\in \L\})
P_{\nur}(\s(\L^c)>k\d) P_{\nur}(\tau>t).
\end{equation}
We establish in the next lemma that $P_{\nur}(\s(\L^c)>k\d)>0$,
which concludes the proof.
\qed

\bl{lem6} Let $\sigma(\L^c)$ be the first time one particle starting
outside $\L$ enters $\L$. Then, for any $\kappa>0$, 
$P_{\nur}(\sigma(\L^c)>\kappa)>0$.
\el
\bpr
We use the coupling described in section~\ref{notations}. Thus, if
$\tilde \sigma(\L^c)$ is the stopping time corresponding to
the coupled independent random walks,
we have $\tilde \sigma(\L^c)\le \sigma(\L^c)$. Thus,
\be{eq3.21}
P_{\nur}(\sigma(\L^c)>\kappa)\ge
P_{\nur}(\tilde \sigma(\L^c)>\kappa)=
\int \prod_{i\not\in \L}\PP(X(i,t)\not\in \L,\ \forall t\le \kappa)^{\eta(i)}
d\nur= \prod_{i\not\in \L}\frac{Z(\gamma  (1-\delta_i))}{Z(\gamma)},
\end{equation}
with $\delta_i=\PP(X(i,t)\in \L,\ t\le \kappa)$.
Now, by Jensen's inequality
\[
\frac{Z(\gamma  (1-\delta ))}{Z(\gamma)}\ge (1-\delta)^{\rho}.
\]
Thus,
\be{eq3.22}
P_{\nur}(\sigma(\L^c)>\kappa)\ge \left(\prod_{i\not\in\L}(1-\delta_i)\right)^
{\rho}>0 \Longleftrightarrow 
\sum_{i\in \ZZ^d} \delta_i<\infty.
\end{equation}
Now, a particle starting on $i$ reaches $\L$ within time $\kappa$,
if it makes at least
$d(i,\L)/R$ jumps within time $\kappa$ (recall that $R$ is the range of
$p$). Thus, if $d(i)$ is the integer part
of $d(i,\L)/R$,
\be{eq3.23}
\PP(X(i,t)\in \L,\ t\le \kappa)\le
\sum_{n\ge d(i)} e^{-\D \kappa}\frac{(\D \kappa)^n}{n!}
\le \frac{(\D \kappa)^{d(i)}}{d(i)!}.
\end{equation}
Hence, the series in (\ref{eq3.22}) is converging.
\epr

\noindent
{\bf Proof of Proposition~\ref{lem5}.}
The proof follows the same arguments as in 
the proof of Theorem 3 c), of \cite{ad},
once the inequality (\ref{eq3.7}) is established with 
$\e_i = \PP_i(S_n \notin \L, \, \forall n \in \NN)$. It goes as follows.
Let $\nu_{\epsilon}$ be the product measure
\[
d \nu_{\epsilon}(\eta) = \prod_{i \in \L} d\theta_{\gamma(\rho)}(\eta_i)
\prod_{i \notin \L} d\theta_{\epsilon_i \gamma(\rho) }(\eta_i)\, .
\]
Let  $\L_n := [-n;n]^d$ and ${\cal{G}}_n$ be the $\sigma$-algebra 
$\sigma(\eta_i; i \in \L_n)$,
then 
\be{densite.nue}
\nur \text{ p.s.} \,\,
\left. \frac{d\nu_{\epsilon}}{d\nur} \right|_{{\cal G}_n}
= \prod_{i \in \L^c\cap\L_n} 
\frac{\epsilon_i^{\eta_i} Z(\gamma)}{Z(\epsilon_i \gamma)} \, ;
\quad
\nu_{\epsilon}\text{ p.s.} \,\, 
\left. \frac{d\nur}{d\nu_{\epsilon}} \right|_{{\cal G}_n}
= \prod_{i \in \L^c\cap\L_n} 
\frac{\epsilon_i^{- \eta_i} Z(\epsilon_i \gamma)}{Z( \gamma)} \, .
\end{equation}
Let $h(\alpha)$ denote the 
Laplace transform of $\theta_{\gamma}$; i.e. 
$h(\alpha) = Z(e^{\alpha}\gamma)/Z(\gamma)$. Note that $h$ is defined for 
any $\alpha$ such that $e^{\alpha}\gamma < \sup g(k)$, and is analytic in
this domain. In particular, $h$ is analytic in a neighbourhood of $0$.
For all $i \notin \L$, let $\alpha_i$ be defined by  
$e^{-\alpha_i}=\epsilon_i$.  A simple computation then yields for 
all $p \ge 1$,
\be{a.c.nue}
\begin{array}{ll}
& \int ( \left.
\frac{d\nu_{\epsilon}}{d\nur} \right|_{{\cal G}_n})^p \, d\nur
= 
\prod_{i \in \L^c \cap \L_n} 
\frac{Z(\epsilon_i^p \gamma)}{Z(\gamma)}
\frac{Z(\gamma)^p}{Z(\epsilon_i \gamma)^p} 
=
\prod_{i \in \L^c \cap \L_n} 
\frac{h(-p\alpha_i)}{h(-\alpha_i)^p} \, ; \\
\text{ and }
&
\int  ( \left. 
\frac{d\nur}{d\nu_{\epsilon}} \right|_{{\cal G}_n} )^p \, d\nu_{\epsilon}
=  
\prod_{i \in \L^c \cap \L_n}  
\frac{Z(\epsilon_i^{-(p-1)} \gamma)}{Z(\gamma)}  
\frac{Z(\epsilon_i \gamma)^{p-1}}{Z(\gamma)^{p-1}} 
= 
\prod_{i \in \L^c \cap \L_n}  
h(\alpha_i(p-1)) h(-\alpha_i)^{p-1} \, . 
\end{array}  
\end{equation}

The functions $m_p: \alpha \mapsto \frac{h(-p\alpha)}{h(-\alpha)^p}$
and $n_p: \alpha \mapsto h(\alpha(p-1)) h(-\alpha)^{p-1}$ are analytic in a 
neighbourhood  of $0$, and satisfy $m_p(0)=n_p(0)=1$, 
$m'_p(0)=n'_p(0)=0$,  $m''_p(0)=n''_p(0) > 0$ for $p>1$. Therefore, 
the products in (\ref{a.c.nue}) have finite limits when $n \rightarrow \infty$,
as soon as $\sum_{i \in \L^c} (1-\epsilon_i)^2 < +\infty$. 
In the asymmetric case, the Fourier transform of the Green function
has a singularity at 0 which is square integrable as soon as $d\ge 3$,
so that the above series is convergent. Thus, for $d \ge 3$, 
$\left .\frac{d\nu_{\epsilon}}{d\nur} \right|_{{\cal G}_n}$
is a $(P_{\nur},({\cal G}_n))$ martingale, which is 
uniformly bounded in $L^p(\nur)$ for all
$p \ge 1$. It  follows from the martingale convergence theorem that 
$\nu_{\epsilon}$  is a.c. with respect to $\nur$, with $\frac{d\nu_{\epsilon}}
{d\nur} \in L^p(\nur)$.  In the same way, $\nur$ is a.c. with respect to 
$\nu_{\epsilon}$, and  $\frac{d\nur}{d\nu_{\epsilon}}\in L^p(\nu_{\epsilon})$.  
 
Following \cite{ad}, we prove that this yields uniform   
$L^p(d\nur)$-estimates of  $f_t :=dT_t(\nur)/ d\nur$, for $p\ge 1$. 
First of all, let us express the density of $\nu_t:=T_t(\nur)$  with respect 
to $\nur$. For $\v$ continuous and bounded
\[
\int \v dT_t(\nur)=\frac{\int\bar S_t(\v) 1_{\A^c}d\nur}
{\int\bar S_t(1_{\A^c}) 1_{\A^c}d\nur}=\int \v
\frac{\bar S_t^*(1_{\A^c})}{P_{\nur}^*(\tau>t)}d\nur,
\]
so that $\nur$-a.s.
$f_t =\frac{P^*_{\eta}(\tau>t)}{P_{\nur}^*(\tau>t)}$. 

Let $\A_0=\{\eta; \, \forall i \in \L, \, \eta_i=0 \}$. We 
prove now that for any increasing function $\v$,
\be{comp.nut.nue}
\int_{\A_0} \v  \, d\nu_t \ge \frac{\nu_t(\A_0)}{\nur(\A_0)} 
\int_{\A_0} \v  \, d\nu_{\epsilon} \, . 
\end{equation} 
To this end, let us write $\eta=(\eta_{\L}, \eta_{\L^c})$ the decomposition
of $\NN^{\ZZ^d}$ in $\NN^{\L} \times \NN^{\L^c}$. Moreover, if 
$\mu$ is a probability measure on $\NN^{\ZZ^d}$, let $\pi_{\L^c}(\mu)$
denote its projection on $\sigma(\eta_i, i\in \L^c)$. We have 
\[
\int_{\A_0} \v  \, d\nu_t
= \nur(\A_0) \int \v(0,\eta_{\L^c}) f_t(0,\eta_{\L^c}) 
	\frac{d\nur}{d\nu_{\epsilon}}(\eta_{\L^c}) 
	\, d\pi_{\L^c}(\nu_{\epsilon}) \, .
\]
By (\ref{eq3.7}), $\forall i \notin \L$, $\AA_i  f_t(0,\eta_{\L^c})
\ge \epsilon_i f_t(0,\eta_{\L^c})$, and $\AA_i 	\frac{d\nur}{d\nu_{\epsilon}}
= \frac{1}{\epsilon_i} \frac{d\nur}{d\nu_{\epsilon}}$. Therefore,
$f_t(0,\eta_{\L^c}) \frac{d\nur}{d\nu_{\epsilon}}(\eta_{\L^c})$ is an 
increasing function of $\eta_{\L^c}$. $\pi_{\L^c}(\nu_{\epsilon})$ being
a product measure, it follows from FKG's inequality that
\[
\int_{\A_0} \v  \, d\nu_t
\ge \nur(\A_0) \int \v(0,\eta_{\L^c}) \, d\pi_{\L^c}(\nu_{\epsilon})
\int f_t(0,\eta_{\L^c}) 
	\frac{d\nur}{d\nu_{\epsilon}}(\eta_{\L^c}) 
	\, d\pi_{\L^c}(\nu_{\epsilon}) \, ,
\]
which is just (\ref{comp.nut.nue}).

We apply now (\ref{comp.nut.nue}) to the decreasing function
$f_t^{p-1} (\frac{d\nu_{\epsilon}}{d\nur})^r$ ($p\ge 1, r\ge0$). We obtain
\[ 
\begin{array}{ll}
\int_{\A_0} f_t^p \left(\frac{d\nu_{\epsilon}}{d\nur} \right)^r \, d\nur
& = \int_{\A_0} f_t^{p-1} 
	\left(\frac{d\nu_{\epsilon}}{d\nur} \right)^r \, d\nu_t \\
& \le \frac{\nu_t(\A_0)}{\nur(\A_0)} 
 \int_{\A_0} f_t^{p-1}  (\frac{d\nu_{\epsilon}}{d\nur})^r \, d\nu_{\epsilon} \\
& \le \frac{\nu_t(\A_0)}{\nur(\A_0)} 
\int_{\A_0} f_t^{p-1}  (\frac{d\nu_{\epsilon}}{d\nur})^{r+1} \, d\nur
\end{array}
\]
It follows by induction that $\forall p, r \ge 0$, 
\[
\int_{\A_0} f_t^p \left(\frac{d\nu_{\epsilon}}{d\nur} \right)^r \, d\nur
\leq \left(\frac{\nu_t(\A_0)}{\nur(\A_0)} \right)^p 
\int_{\A_0} \left(\frac{d\nu_{\epsilon}}{d\nur}\right)^{p+r}  \, d\nur \, .
\]
Taking $r=0$, and applying once more FKG's inequality to the decreasing
functions $1_{{\A_0}}$ and  $f_t^p$, we get $\forall p \ge 1$,
\[
\nur(\A_0) \int  f_t^p \, d\nur
\leq \int_{\A_0} f_t^p \, d\nur 
\leq \left(\frac{\nu_t(\A_0)}{\nur(\A_0)} \right)^p 
\int_{\A_0} \left(\frac{d\nu_{\epsilon}}{d\nur} \right)^p \, d\nur \, ,
\]
so that $\forall p \geq 1$, 
\be{ft-uLp}
\sup_t \int  f_t^p \, d\nur \leq  \frac{1}{\nur(\A_0)^{p+1}} 
\int_{\A_0} \left(\frac{d\nu_{\epsilon}}{d\nur}\right)^p \, d\nur \, .
\end{equation}

This in turn implies uniform $L^p(\nur)$-estimates for 
$\frac{d\Phi^n(\nur)}{d\nur}$. Indeed, using expression (\ref{eq.nun}), 
if we define
\be{eq3.24}
dm_n(t)=\frac{\Pn t^ndt}{\int_0^{\infty} \Pn t^ndt},\quad
\text{ then}\quad \frac{d\Phi^n(\nur)}{d\nur}=
\int_0^{\infty} \frac{dT_t(\nur)}{d\nur}dm_{n-1}(t).
\end{equation}
Thus, using H\"older inequality for $p\ge 1$, 
\be{eq3.25}
\sup_{t>0} \int\left(\frac{dT_t(\nur)}{d\nur}\right)^pd\nur\le C
\Longrightarrow 
\sup_{n} \int\left(\frac{d\Phi^n(\nur)}{d\nur}\right)^pd\nur\le C.
\end{equation}
Moreover, we obtain the same uniform bounds for the Cesaro limit,
and Proposition~\ref{lem5} follows.
\qed

\noindent
{\bf Proof of Corollary~\ref{cor1}.}
We define the map $\Phi_*$ associated to the time reversed dynamics.
If $\nu$ is such that $E_{\nu}^*[\tau]<\infty$, then
\[
\int \v d\Phi_*(\nu)=\frac{1}{E_{\nu}^*[\tau]}\int_0^{\infty}
\int \bar S_t^*(\v) d\nu dt.
\]
Our previous results (Proposition~\ref{lem5}) hold equally for
$\bar \nu_n^*:=\frac{1}{n}\left( \Phi_*(\nur)+\dots+\Phi_*^n(\nur)\right)$,
with the consequences that $\{\bar \nu_n^*, n\in \NN\}$ is tight and
$g_n:=d\bar \nu_n^*/d\nur$ is uniformely in $L^p(\nur)$ for any
$p\ge 1$ in dimensions
$d\ge 3$. Let $f_n$ be the density of $\bar \nu_n$ with respect to $\nur$,
and assume that $\{f_n\}$ converge along a subsequence $\{n_k\}$ to $f$ 
solution of (\ref{eq.eigen}) and that $\{g_n\}$ converge along a subsequence 
$\{m_i\}$ to $g$ solution to the adjoint equation to (\ref{eq.eigen}).
We can as well assume that these convergence hold in weak $L^2(\nur)$.
As $f_{n}$ and $g_n$ are nonincreasing functions, we have by FKG's inequality
\[
\int f_{n_k} g_{m_i}d\nur \ge \int f_{n_k} d\nur \int g_{m_i}d\nur=1 .
\]
After taking first the limit in $k$, and then in $i$, we obtain
$\int fg d\nur\ge 1$. Also, this integral is finite by Cauchy-Schwarz.
Thus, we can define $d\tilde \nur= fgd\nur/(\int fgd\nur)$. Let
$dQ_t(\eta_.)$ be the probability measure on paths, defined by
\be{eq6.1}
dQ_t(\eta_.):= \frac{e^{\lambda(\rho) t} g(\eta_t) f(\eta_0)}{\int fgd\nur}
1_{\tau>t} dP_{\nur}(\eta_.).
\end{equation}
For $\v$ such that $\v g\in L^2(\nur)$, we obtain using (\ref{eq.eigen}),
\begin{align}
\int \v(\eta_t)dQ_t(\eta_.)&=\frac{\int
E_{\eta}[\v(\eta_t)g(\eta_t)1_{\tau>t}]f(\eta)e^{\lambda(\rho) t} d\nur(\eta)}
{\int fgd\nur}=\frac{\int\bar S_t(\v g) f e^{\lambda(\rho) t} d\nur}
{\int fgd\nur} \notag\\
&=\frac{\int\v g \bar S_t^*(f)e^{\lambda(\rho) t} d\nur}{\int fgd\nur}
=\int \v d\tilde \nur.\notag
\end{align}
Also, if $\v$ is such that $\v f\in L^2(\nur)$,
\[
\int \v(\eta_0)dQ_t(\eta_.)
=\frac{\int \bar S_t(g) \v f e^{\lambda(\rho) t} d\nur}
{\int fgd\nur}=\int \v d\tilde \nur.
\]
Now, by applying Jensen's inequality and recalling 
that $f,g\in L^p(\nur)$ for $p\ge 1$,
\begin{align}
\log(P_{\nur}(\tau>t))&=\log(\int fgd\nur)+
\log\left(\int \frac{e^{-\lambda(\rho) t}}
{g(\eta_t)f(\eta_0)}dQ_t(\eta_.)\right)
\notag\\
&\ge \log(\int fgd\nur)-\int \log\left( g(\eta_t)\right) dQ_t(\eta_.)-
\int\log\left( f(\eta_0)\right) dQ_t(\eta_.)-\lambda(\rho) t\notag\\
&\ge \log(\int fgd\nur)-\int \log(fg) d\tilde \nur -\lambda(\rho) t.\notag
\end{align}
This concludes the proof of the Corollary.
\qed
\section{Example.}
\label{example}
Let us consider the totally asymmetric simple exclusion
in one dimension. Thus,
\[
\forall i\in \ZZ,\quad p(i,i+1)=1,\quad{\rm and}\quad
p(i,j)=0\ \text{if }j\not=i+1.
\]
Let $\tau$ be the first time the origin is occupied.
Let $\chi(\eta):=\inf\{k\ge 0:\ \eta(-k)=1\}$, and $N_t$ be a Poisson process
of intensity 1. A simple computation yields 
\be{eq5.1}
P_{\nur}(\tau>t)=\int\PP(N_t<\chi(\eta))d\nur(\eta)=
\sum_{k=1}^{\infty} \rho(1-\rho)^k \PP(N_t<k)=(1-\rho)e^{-\rho t}.
\end{equation}
Thus,
\be{eq5.2}
\frac{P_{\nur}(\tau>t+s)}{P_{\nur}(\tau>t)}=e^{-\rho s}\quad
\text{and}\quad 
\lambda(\rho):= \lim_t -\frac{1}{t} \log\left(P_{\nur}(\tau>t)\right)=\rho.
\end{equation}
Following the approach of the proof of Theorem 3c) of \cite{ad},
it is easy to establish that the Yaglom limit exists and is
\be{eq5.3}
d\mur(\eta)=\prod_{i<0}d\B_{\rho}(\eta_i)\prod_{i\ge 0} d\B_0(\eta_i)\quad
\text{where $\B_{\rho}$ is the Bernoulli probability of parameter $\rho$}.
\end{equation}
Can we approximate $\mur$ and $\lambda(\rho)$ by the corresponding quantities
for the process on a large circle? The answer is no, as we shall see.

Let $\C_N=\{0,1,\dots,N\}$ where
sites $N$ and 0 are identified, and consider the generator
\be{eq5.4}
\LL_N\v=\sum_{i=0}^{N-1}\eta(i)\left(1-\eta(i+1)\right)
\left(\v(\eta^i_{i+1})-\v(\eta)\right),
\end{equation}
with as invariant measure $\nu_N$, which is the uniform measure 
on all configurations
with $[\rho N]$ particles on $\C_N$. 

Let $P_{\eta,N}$ be the law of the process generated by $\LL_N$, and
let $\eta$ be in the support of $\nu_N$. Then,
\be{eq5.5}
P_{\eta,N}(\tau>t)=e^{-t} 
\sum_{k=1}^{\chi(\eta)-1} \frac{t^k}{k!}.
\end{equation}
Thus, for a polynomial $Q_N$ of degree at most $N$
\be{eq5.6}
P_{\nu_N,N}(\tau>t)=e^{-t}Q_N(t)
\quad\Longrightarrow\quad
\lambda_N(\rho):=\lim_t -\frac{1}{t} \log\left(P_{\nu_N,N}(\tau>t)\right)=1.
\end{equation}
Also, it is an easy computation which yields 
\be{eq5.7}
\lim_t \frac{P_{\eta,N}^*(\tau>t)}{P_{\nu_N,N}^*(\tau>t)}=
\left( \begin{array}{c}
N \\
\left[ \rho N \right]
\end{array} \right)
\prod_{i=1}^{[\rho N]} \eta(-i)
\quad\text{ and }\quad
\lim_t \frac{P_{\nu_N,N}(\tau>t+s)}{P_{\nu_N,N}(\tau>t)}=e^{-s}.
\end{equation}
Thus, as in \cite{ad}, one concludes the existence
of a Yaglom limit $\mu_N$ concentrated on the configurations with
particles occupying all $[\rho N]$ sites to the ``left'' of 0.
Thus, $\mu_N$ and $\lambda_N(\rho)$ do converge, 
but to $\mu_1$ and 1 respectively, and this approach
misses all the $\mu_{\rho}$ with $\rho<1$.

\end{document}